\documentclass[11pt]{article}
\usepackage[a4paper,margin=1in]{geometry}
\usepackage{amsmath,amssymb,amsthm,mathtools}
\usepackage{newtxtext,newtxmath}
\usepackage{xcolor}
\definecolor{paperlinkblue}{HTML}{174A76}
\usepackage[unicode,colorlinks=true,bookmarksnumbered,bookmarksopen=true,linktoc=all]{hyperref}
\hypersetup{
  linkcolor=paperlinkblue,
  citecolor=paperlinkblue,
  urlcolor=paperlinkblue,
  pdftitle={Poincar\'e-Type Hessian Quotient Equations},
  pdfauthor={Linwei Zhang}
}
\numberwithin{equation}{section}
\newtheorem{theorem}{Theorem}[section]
\newtheorem{lemma}[theorem]{Lemma}
\newtheorem{proposition}[theorem]{Proposition}
\theoremstyle{definition}
\newtheorem{definition}[theorem]{Definition}
\newtheorem{remark}[theorem]{Remark}
\newcommand{\ii}{\sqrt{-1}}
\newcommand{\ddc}{dd^c}

\newcommand{\PM}{\mathrm{PM}}

\newcommand{\Vol}{\operatorname{Vol}}

\newcommand{\ee}{\varepsilon}
\title{Poincar\'e-Type Hessian Quotient Equations}
\author{Linwei Zhang}
\date{16 October 2026}
\begin{document}
\maketitle
\begin{abstract}
Inspired by the two-parameter continuity-path approach of Chen and
Xu~\cite{ChenXu}, we employ a modified continuity path to characterize
the solvability of the Hessian quotient equation for K\"ahler metrics
with Poincar\'e-type singularities along a divisor $D$, allowing simple
normal crossings and self-intersections.
\end{abstract}
\tableofcontents
\clearpage
\section{Introduction}

The search for canonical K\"ahler metrics has long been a central theme
in K\"ahler geometry. In his celebrated resolution of the Calabi
conjecture, Yau~\cite{YauCalabi} proved the existence of K\"ahler metrics
with prescribed Ricci form on compact K\"ahler manifolds. In particular,
together with Aubin's independent work in the negative first Chern class
case, this established the existence of K\"ahler--Einstein (KE) metrics
when the first Chern class is negative or vanishes. Extending the
Ricci-flat theory beyond the compact setting, Tian and Yau
\cite{TianYauI,TianYauII} constructed complete Ricci-flat K\"ahler
metrics on suitable noncompact K\"ahler manifolds in their works
published in 1990 and 1991. Their methods laid important analytic
groundwork for the study of canonical metrics on complete manifolds.

Beginning with his work in 2011, Auvray initiated a systematic study of
K\"ahler metrics with Poincar\'e-type cusp singularities, including
topological constraints for constant scalar curvature K\"ahler metrics,
the geometry of the space of Poincar\'e-type metrics, and the asymptotic
behavior of extremal K\"ahler metrics~\cite{AuvrayTopological,Auvray,AuvrayExtremal}.
This direction has subsequently attracted considerable attention. More recently, Chen and
Xu~\cite{ChenXu} introduced the $J$-equation into this cusp setting,
building on the earlier compact theory of Song and
Weinkove~\cite{SongWeinkove}. They showed that the subsolution method
continues to play a decisive role: solvability of the Poincar\'e-type
$J$-equation is characterized by the existence of a subsolution
satisfying both a strict cone condition and an appropriate asymptotic
condition near the divisor. Their argument draws on the analytic
framework for complete manifolds developed in the work of Tian--Yau
and Auvray.

Motivated by the recent work of Chen, Nie and Xu~\cite{ChenNieXu}
on complex Hessian-type equations, as well as the results of
Chen--Xu~\cite{ChenXu} for the $J$-equation, we aim to establish an
analogous solvability criterion for Poincar\'e-type Hessian quotient
equations in the cusp setting. Specifically, we seek to characterize
solvability by the existence of a subsolution satisfying an endpoint
cone condition and a prescribed asymptotic condition. The equation
and the precise hypotheses are stated below.

\section{Main theorem}
\label{sec:intro}
Let $X$ be a connected compact K\"ahler manifold of complex dimension
$n\ge2$, let $D\subset X$ be a simple normal crossing divisor (allowing
normal self-intersections), and put $M=X\setminus D$.
Fix a Poincar\'e-type K\"ahler class $\Omega$, a representative
$\omega\in\PM_\Omega$, and a Poincar\'e-type K\"ahler form $\chi$.
Throughout, $1\le k<n$ and $m=n-k$.
We use the volume normalization
\begin{equation}\label{eq:intro-normalization}
  \int_M\omega^n=\int_M\omega^k\wedge\chi^m=1.
\end{equation}
 We consider
\begin{equation}\label{eq:intro-target}
  \omega_\varphi^n=\omega_\varphi^k\wedge\chi^m,
 \qquad \omega_\varphi=\omega+\ddc\varphi>0.
\end{equation}
Solutions are understood as Poincar\'e-type metrics in $\PM_\Omega$;
their potentials are unique only up to additive constants.

Our solvability conditions require only the \emph{endpoint} Hessian
quotient. More precisely, there must exist
$\omega\in\PM_\Omega$ and constants $\delta,\eta_0,C>0$ such that
\begin{align}
 n\omega^{n-1}&\ge(1+\delta)k\,\omega^{k-1}\wedge\chi^m,
 &&\tag{B}\label{eq:intro-B}\\
 \left|\frac{\omega^k\wedge\chi^m}{\omega^n}-1\right|
 &\le C\rho^{-\eta_0},
 &&\tag{C}\label{eq:intro-C}
\end{align}
where $\rho$ is the cusp weight of Section~\ref{sec:preliminaries}.
The inequality in (B) is in the sense of strongly positive
$(n-1,n-1)$-forms. Neither condition is imposed for any lower
mixed power $\chi^j$, $1\le j<m$.

\begin{theorem}[Solvability criterion for Poincar\'e-type Hessian quotient equations]
\label{thm:main}
The equation \eqref{eq:intro-target} admits a Poincar\'e-type
K\"ahler solution in $\PM_\Omega$ \emph{if and only if} there exists
a representative $\omega\in\PM_\Omega$ satisfying \emph{(B)} and
\emph{(C)}. The solution metric is unique in $\PM_\Omega$ and its
potential is unique modulo constants.
\end{theorem}

The necessity of these endpoint conditions is elementary. Indeed, if
$g\in\PM_\Omega$ solves \eqref{eq:intro-target}, choose $\omega=g$.
Then (C) holds with zero error. To verify (B), diagonalize $\chi$
relative to $g$, with eigenvalues $\beta_i>0$. The equation says
$\sigma_m(\beta)=\binom nm$. On the hyperplane perpendicular to the
$i$th eigenvector, the ratio of $k g^{k-1}\wedge\chi^m$ to
$n g^{n-1}$ equals
\[
 \frac{\sigma_m(\beta|i)}{\binom nm}
 =1-\frac{\beta_i\sigma_{m-1}(\beta|i)}{\binom nm}<1.
\]
Since $g$ and $\chi$ are uniformly equivalent Poincar\'e metrics,
this inequality has a uniform strict margin, yielding (B).
The sufficiency direction is approached by the continuation path below.

\phantomsection
\addcontentsline{toc}{subsection}{The Chen--Xu two-parameter method and our path}
\subsection*{The Chen--Xu two-parameter method and our path}
The two-parameter method of Chen--Xu~\cite{ChenXu} introduces a small
positive $\varepsilon$ to make the potential linearization invertible.
For general Hessian quotient equations, however, interpolating the
$(1,1)$-forms $\omega$ and $\chi$ before taking the $m$th wedge power
leads to unnecessary intermediate mixed cone and decay assumptions.
Instead, we interpolate their closed positive $(m,m)$-powers directly:
\begin{equation}\label{eq:intro-theta}
  \Theta_t=(1-t)\omega^m+t\chi^m,
 \qquad 0\le t\le1.
\end{equation}
Our perturbed continuity path is
\begin{equation}\label{eq:intro-path}
  \omega_{\varphi_{\varepsilon,t}}^n
 =e^{\varepsilon\varphi_{\varepsilon,t}}
 \omega_{\varphi_{\varepsilon,t}}^k\wedge\Theta_t,
 \qquad (\varepsilon,t)\in(0,\varepsilon_*)\times[0,1].
\end{equation}
At $t=0$, $\varphi_{\varepsilon,0}=0$ is a solution; at $t=1$
we obtain the $\varepsilon$-perturbation of
\eqref{eq:intro-target}. By the normalization
\eqref{eq:intro-normalization} and Gaffney--Stokes,
\[
 \int_M\omega_\varphi^k\wedge\Theta_t
 =(1-t)\int_M\omega^n+t\int_M\omega^k\wedge\chi^m=1,
\]
so no path-dependent normalizing coefficient is needed. Moreover,
(B) implies, for some $\delta_*>0$ independent of $t$,
\[
 n\omega^{n-1}\ge(1+\delta_*)k\omega^{k-1}\wedge\Theta_t,
\]
whereas (C) gives
\[
 \left|\frac{\omega^k\wedge\Theta_t}{\omega^n}-1\right|
 =t\left|\frac{\omega^k\wedge\chi^m}{\omega^n}-1\right|
 \le C\rho^{-\eta_0}.
\]
The operator along this path is a positive mixture of two inverse
Hessian quotients. Sections~\ref{sec:openness} and~\ref{sec:C0}
prove openness and the uniform $C^0$ estimates for this new path.
Higher-order estimates and final continuation require the corresponding
adaptation of Sections~\ref{sec:high-order} and~\ref{sec:main-proof}.
\\
\phantomsection
\addcontentsline{toc}{subsection}{Acknowledgements}
\paragraph{Acknowledgements.}
The author would like to express his sincere gratitude to Professor Haozhao Li
for his longstanding support and encouragement. The author is also grateful
to Yulun Xu, Sijie Nie, and Yongheng Han for many helpful discussions.

\phantomsection
\addcontentsline{toc}{subsection}{On the Use of AI}
\paragraph{On the Use of AI.}
The mathematical ideas and research direction of this paper are the author's
own. In particular, the modified continuity path
\[
\Theta_t=(1-t)\omega^{n-k}+t\chi^{n-k},
\]
in contrast to the path based on $\chi_t=(1-t)\omega+t\chi$, was developed
by the author through independent experimentation. The author used ChatGPT
Plus to explore an alternative to Lemma~4.7 of Chen--Xu~\cite{ChenXu}
and to assist in adapting the proof of Proposition~4.5 in that work to the
modified path. In particular, ChatGPT Plus played a key role in working
out the arguments in Steps~2 and~3 of Proposition~\ref{prop:C2} of the
present paper. The author has reviewed the mathematical arguments and
takes full responsibility for the contents of the paper and for any
remaining errors.

\section{Poincar\'e geometry and analytic preliminaries}
\label{sec:preliminaries}
Let $(X,\omega_X)$ be a compact K\"ahler manifold of complex dimension $n$,
let $D\subset X$ be a divisor with simple normal crossings, possibly including
normal self-intersections, and put $M=X\setminus D$.
We use the convention $\ddc=\ii\partial\bar\partial$.
We follow the notation of Chen--Xu~\cite{ChenXu}.

\subsection{Poincar\'e-type K\"ahler metrics and the space \texorpdfstring{$\PM_\Omega$}{PM Omega}}
\label{subsec:PM}
Around any $p\in D$, choose a cusp coordinate chart $(U;z_1,\ldots,z_n)$
such that $D\cap U=\{z_1\cdots z_r=0\}$. The model Poincar\'e metric is
\begin{equation}\label{eq:poincare-model}
 \omega_{\mathrm P}=\ii\sum_{j=1}^{r}
 \frac{dz_j\wedge d\bar z_j}{|z_j|^2(\log|z_j|^2)^2}
 +\ii\sum_{j=r+1}^{n}dz_j\wedge d\bar z_j.
\end{equation}

\begin{definition}[Poincar\'e-type K\"ahler metric]\label{def:PM}
Let $\Omega=[\omega_X]\in H^2(X,\mathbb R)$. A K\"ahler form
$\omega$ on $M$ belongs to $\PM_\Omega$ if
\begin{enumerate}
 \item $\omega$ is smooth on $M$ and, in each cusp chart,
 $C^{-1}\omega_{\mathrm P}\le \omega\le C\omega_{\mathrm P}$;
 \item $\omega=\omega_X+\ddc\varphi$ for some $\varphi\in C^\infty(M)$ satisfying
 \[
 \varphi=O\!\left(1+\sum_{j=1}^{r}\log(-\log|z_j|^2)\right),\qquad
 |\nabla_{\omega_{\mathrm P}}^\ell\varphi|_{\omega_{\mathrm P}}\le C_\ell,
 \quad \ell\ge1.
 \]
\end{enumerate}
\end{definition}
Every $\omega\in\PM_\Omega$ is complete and has finite volume.
Fix $\omega\in\PM_\Omega$. We denote the corresponding space of
Poincar\'e-type potentials by
\[
 \widetilde{\PM}_\Omega
 =\{\varphi\in C^\infty(M,\mathbb R):\omega+\ddc\varphi\in\PM_\Omega\}.
\]

\subsection{Quasi-coordinates}
\label{subsec:qc}
Let $\Delta$ be the unit disc and $\Delta^*=\Delta\setminus\{0\}$.
For $0<\delta<1$, define
\begin{equation}\label{eq:quasi-map}
 \phi_\delta:\tfrac34\Delta\longrightarrow\Delta^*,\qquad
 \phi_\delta(\xi)=\exp\!\left(-\frac{1+\delta}{1-\delta}
 \frac{1+\xi}{1-\xi}\right).
\end{equation}
For $\boldsymbol\delta=(\delta_1,\ldots,\delta_r)\in(0,1)^r$,
the associated product quasi-coordinate map is
\begin{equation}\label{eq:Phi}
 \Phi_{\boldsymbol\delta}:P:=(\tfrac34\Delta)^r\times\Delta^{n-r}
 \longrightarrow(\Delta^*)^r\times\Delta^{n-r},
 \quad (\xi_1,\ldots,\xi_n)\longmapsto
 (\phi_{\delta_1}(\xi_1),\ldots,\phi_{\delta_r}(\xi_r),
 \xi_{r+1},\ldots,\xi_n).
\end{equation}
The pullbacks $\Phi_{\boldsymbol\delta}^*\omega_{\mathrm P}$ are
uniformly equivalent to the Euclidean metric on fixed relatively compact
sub-polydiscs, independently of $\boldsymbol\delta$.
For $\omega\in\PM_\Omega$ the same uniform equivalence holds for
$\Phi_{\boldsymbol\delta}^*\omega$.

\subsection{The weight \texorpdfstring{$\rho$}{rho} and function spaces}
\label{subsec:spaces}
Choose finitely many cusp charts $U_i$ covering $D$, with
\[
 D\cap U_i=\{z_{i,1}\cdots z_{i,r_i}=0\},\qquad 1\le i\le N,
\]
and an additional relatively compact interior chart $U_0\Subset M$.
Let $\{\chi_0,\chi_1,\ldots,\chi_N\}$ be a smooth partition of unity
subordinate to this cover. Shrinking the cusp charts if necessary, assume
$-\log|z_{i,j}|^2\ge1$ on their supports. Set
\begin{equation}\label{eq:rho}
 \rho_i=\prod_{j=1}^{r_i}(-\log|z_{i,j}|^2),\qquad
 \rho=\chi_0+\sum_{i=1}^{N}\chi_i\rho_i. 
\end{equation}
Thus $\rho\ge1$ and, on each sufficiently small cusp chart,
\[
 \rho\approx\prod_{j=1}^{r}(-\log|z_j|^2),\qquad
 \Phi_{\boldsymbol\delta}^*\rho
 \approx\prod_{j=1}^{r}(1-\delta_j)^{-1}.
\]
The constants in $\asymp$ are uniform in $\boldsymbol\delta$.

\begin{definition}[Quasi-coordinate H\"older spaces]\label{def:holder}
For $m\in\mathbb N_0$ and $0\le\alpha<1$, set
\[
 \|f\|_{C^{m,\alpha}(U_i\setminus D)}
 :=\sup_{\boldsymbol\delta\in(0,1)^{r_i}}
 \|\Phi_{i,\boldsymbol\delta}^*f\|_{C^{m,\alpha}(P)}.
\]
With the above partition of unity, define
\begin{equation}\label{eq:holder}
 \|f\|_{C^{m,\alpha}(M)}
 :=\|\chi_0 f\|_{C^{m,\alpha}(U_0)}
 +\max_{1\le i\le N}
 \|\chi_i f\|_{C^{m,\alpha}(U_i\setminus D)}.
\end{equation}
\end{definition}

\begin{definition}[Weighted H\"older spaces]\label{def:weighted-holder}
For $\eta\in\mathbb R$, define
\begin{equation}\label{eq:weighted-holder}
 \|f\|_{C^{m,\alpha}_\eta(M)}
 :=\|\chi_0f\|_{C^{m,\alpha}(U_0)}
 +\max_{1\le i\le N}\sup_{\boldsymbol\delta\in(0,1)^{r_i}}
 \left\|\left(\prod_{j=1}^{r_i}(1-\delta_j)^\eta\right)
 \Phi_{i,\boldsymbol\delta}^*(\chi_i f)\right\|_{C^{m,\alpha}(P)}.
\end{equation}
Equivalently, the cusp term may be written using
$\Phi_{i,\boldsymbol\delta}^*(\rho^{-\eta}\chi_i f)$.
We write
\[
 C^\infty_\eta(M)=\bigcap_{m\ge0}C^{m,\alpha}_\eta(M),\qquad
 C^\infty(M)=\bigcap_{m\ge0}C^{m,\alpha}(M),
\]
where a fixed $0<\alpha<1$ is understood.
\end{definition}

\begin{definition}[Weighted Sobolev spaces]\label{def:sobolev}
For $m\in\mathbb N_0$, $1\le p<\infty$, and $\eta\in\mathbb R$, set
\begin{equation}\label{eq:weighted-sobolev}
 W^{m,p}_\eta(M,\omega)
 :=\left\{f\in W^{m,p}_{\mathrm{loc}}(M):
 \|f\|_{W^{m,p}_\eta}^p
 :=\sum_{\ell=0}^{m}\int_M
 |\nabla_\omega^\ell f|_\omega^p\rho^{-\eta p}\,\omega^n<\infty
 \right\}.
\end{equation}
In particular, $W^{m,p}_\eta\subset W^{m,p}_{\eta'}$ if $\eta\le\eta'$.
\end{definition}

\subsection{Poincar\'e and weighted Sobolev inequalities}
\label{subsec:ineq}
We record the following inequalities from Auvray~\cite{Auvray},
as used in Chen--Xu~\cite{ChenXu}.

\begin{lemma}[Poincar\'e inequality]\label{lem:poincare}
There exists $C_P>0$ such that, for every $f\in W^{1,2}_0(M,\omega)$,
\begin{equation}\label{eq:poincare}
 \int_M|f-\bar f|^2\omega^n
 \le C_P\int_M|df|_\omega^2\omega^n,\qquad
 \bar f=\frac{1}{\Vol_\omega(M)}\int_M f\,\omega^n.
\end{equation}
\end{lemma}

\begin{lemma}[Weighted Sobolev inequality]\label{lem:weighted-sob}
Let $m\ge0$ and $1<p,q<\infty$ satisfy
$\frac1q<\frac1p+\frac1{2n}$. Then
\begin{equation}\label{eq:weighted-embedding}
 \|f\|_{W^{m,p}_{-1/p}(M,\omega)}
 \le C\|f\|_{W^{m+1,q}_{-1/q}(M,\omega)}.
\end{equation}
\end{lemma}

\subsection{Stokes' theorem and a generalized maximum principle}
\label{subsec:stokes-max}

\begin{theorem}[Gaffney--Stokes theorem~\cite{Gaffney}]\label{thm:gaffney}
Let $(N,g)$ be a complete oriented Riemannian manifold of real dimension
$d$. If a $C^1$ $(d-1)$-form $\Theta$ satisfies
$|\Theta|_g,|d\Theta|_g\in L^1(N,g)$, then
\begin{equation}\label{eq:gaffney}
 \int_N d\Theta=0.
\end{equation}
\end{theorem}

In particular, Gaffney--Stokes applies to the standard integration-by-parts
forms associated with Poincar\'e-type potentials. If $\omega_u,\omega,\chi$
are Poincar\'e-type K\"ahler metrics and $\omega_u=\omega+\ddc u$, then
\begin{equation}\label{eq:mixed-volume}
 \int_M\omega_u^a\wedge\chi^{n-a}
 =\int_M\omega^a\wedge\chi^{n-a},\qquad 0\le a\le n.
\end{equation}

The exhaustion function $\tau=\log\rho$ is proper on $M$ and has bounded
complex Hessian with respect to the fixed Poincar\'e geometry.
\begin{proposition}[Generalized maximum principle]\label{prop:gen-max}
Let $L=a^{i\bar j}\partial_i\partial_{\bar j}$ be a linear elliptic
operator whose coefficients are uniformly bounded relative to a fixed
Poincar\'e metric $\omega$. If $h\in C^2(M)$ is bounded above, then there
exist $x_\ell\in M$ such that
\begin{equation}\label{eq:max-sequence}
 h(x_\ell)\longrightarrow\sup_Mh,\qquad
 \limsup_{\ell\to\infty}Lh(x_\ell)\le0.
\end{equation}
The analogous assertion holds at the infimum of a bounded-below function.
Consequently, if $b\ge b_0>0$ and $h$ is bounded with
$(L-b)h=f$, then
\begin{equation}\label{eq:max-estimate}
 \|h\|_{L^\infty(M)}\le b_0^{-1}\|f\|_{L^\infty(M)}.
\end{equation}
\end{proposition}

\section{Openness of the continuity path}
\label{sec:openness}
Let $\omega\in\PM_\Omega$ and $\chi$ be the fixed Poincar\'e-type
K\"ahler forms satisfying \eqref{eq:intro-normalization}. Put
$m=n-k$ and
\begin{equation}\label{eq:chi-t}
 \Theta_t=(1-t)\omega^m+t\chi^m,
 \qquad 0\le t\le1.
\end{equation}
The $\Theta_t$ are uniformly positive, closed $(m,m)$-forms and
have uniformly bounded coefficients in quasi-coordinates.
For fixed $\ee>0$ consider
\begin{equation}\label{eq:continuity-path}
 \omega_{\varphi_t}^n
 =e^{\ee\varphi_t}\omega_{\varphi_t}^k\wedge\Theta_t,
 \qquad\omega_{\varphi_t}=\omega+\ddc\varphi_t>0.
\end{equation}
At $t=0$, $\varphi_0=0$ solves the equation. Define
\[
 \mathcal T_{\ee}=
 \{t\in[0,1]:\ \eqref{eq:continuity-path}
 \text{ admits a Poincar\'e-type solution}\}.
\]

\begin{proposition}[Openness]\label{prop:openness}
For every fixed $\ee>0$, $\mathcal T_{\ee}$ is relatively open
in $[0,1]$. At any solution $\varphi_{t_0}$, for each integer
$\ell\ge0$ and $0<\alpha<1$, the potential linearization is
an isomorphism
\begin{equation}\label{eq:linear-isomorphism}
 \mathcal L_{t_0}:C^{\ell+2,\alpha}(M)
 \longrightarrow C^{\ell,\alpha}(M)
\end{equation}
in the unweighted quasi-coordinate H\"older spaces of
Section~\ref{subsec:spaces}.
\end{proposition}

\begin{proof}
Write $g=\omega_{\varphi_{t_0}}$ and define
\begin{equation}\label{eq:log-quotient}
 \mathcal F(t,u)=\log\frac{\omega_u^n}
 {\omega_u^k\wedge\Theta_t}-\ee u.
\end{equation}
The equation is $\mathcal F(t,u)=0$. Since $\Theta_t$ is independent
of $u$, differentiation in a potential direction $v$ gives
\begin{equation}\label{eq:linearization}
 D_u\mathcal F(t_0,\varphi_{t_0})[v]
 =(L_{t_0}-\ee)v,
 \quad L_{t_0}v=\Delta_g v
 -k\frac{\ddc v\wedge g^{k-1}\wedge\Theta_{t_0}}
 {g^k\wedge\Theta_{t_0}}.
\end{equation}
We verify uniform ellipticity even though $\Theta_t$ need not be
an $m$th power. For $\zeta=\omega,\chi$ let
\[
 L_\zeta v=\Delta_g v
 -k\frac{\ddc v\wedge g^{k-1}\wedge\zeta^m}
 {g^k\wedge\zeta^m},
\]
and set
\[
 p_0=\frac{(1-t_0)g^k\wedge\omega^m}
 {g^k\wedge\Theta_{t_0}},\qquad
 p_1=\frac{t_0g^k\wedge\chi^m}
 {g^k\wedge\Theta_{t_0}}.
\]
These are nonnegative smooth functions with $p_0+p_1=1$, and
\begin{equation}\label{eq:ellipticity}
  L_{t_0}=p_0L_\omega+p_1L_\chi.
\end{equation}
For each positive K\"ahler form $\zeta$, the linearization of
$\log(g^n/(g^k\wedge\zeta^m))$ is elliptic. Indeed, in
$\zeta$-unitary coordinates diagonalizing $g$, with eigenvalues
$\lambda_i>0$, its diagonal coefficients are
\[
 F_{\zeta,i}=
 \frac1{\lambda_i}-\frac{\sigma_{k-1}(\lambda|i)}
 {\sigma_k(\lambda)}
 =\frac{\sigma_k(\lambda|i)}
 {\lambda_i\sigma_k(\lambda)}>0.
\]
At the fixed Poincar\'e-type solution, $g,\omega,\chi$ are
uniformly equivalent and have bounded quasi-coordinate derivatives.
Thus both $L_\omega$ and $L_\chi$ are uniformly elliptic with
uniformly controlled coefficients, and the convex combination
\eqref{eq:ellipticity} has the same properties. Hence
$L_{t_0}-\ee$ is uniformly elliptic with strictly negative
zeroth-order term.

The generalized maximum principle (Proposition~\ref{prop:gen-max})
gives, for bounded $v$ satisfying $(L_{t_0}-\ee)v=f$,
\begin{equation}\label{eq:lin-c0}
 \|v\|_{L^\infty(M)}\le\ee^{-1}\|f\|_{L^\infty(M)}.
\end{equation}
In particular the operator is injective. For surjectivity, solve
the zero-boundary Dirichlet problems on a smooth exhaustion of $M$.
The same $L^\infty$ estimate, followed by interior Schauder estimates
on fixed smaller quasi-coordinate polydiscs, gives a limit
$v\in C^{\ell+2,\alpha}(M)$ satisfying
\[
 (L_{t_0}-\ee)v=f,\qquad
 \|v\|_{C^{\ell+2,\alpha}(M)}
 \le C_{\ee}\|f\|_{C^{\ell,\alpha}(M)}.
\]
The constants in these local estimates are uniform across the
quasi-coordinate atlas. Uniqueness follows from
\eqref{eq:lin-c0}, proving \eqref{eq:linear-isomorphism}.

Since $g\approx\omega$ and $\Theta_{t_0}\approx\omega^m$,
the equation itself implies that the fixed solution
$\varphi_{t_0}$ is bounded for $\ee>0$, and hence belongs to the
unweighted quasi-coordinate spaces. The map $\mathcal F$ is smooth
near $(t_0,\varphi_{t_0})$; positivity of $\Theta_t$ persists for
$t$ in a small real interval around $t_0$ if necessary. The Banach
implicit function theorem now supplies nearby Poincar\'e-type
solutions. Thus $\mathcal T_\ee$ is relatively open.
\end{proof}

\begin{remark}\label{rem:openness-eps}
The inverse bound \eqref{eq:lin-c0} depends on $\ee>0$.
At $\ee=0$ constants lie in the kernel of the potential
linearization, so this argument asserts openness only at each fixed
positive $\ee$.
\end{remark}

\section{Uniform \texorpdfstring{$C^0$}{C0} estimates for the mixed path}
\label{sec:C0}
Let $1\le k<n$, $m=n-k$, and assume
\eqref{eq:intro-normalization}. For the path
\begin{equation}\label{eq:C0-normalized-path}
 \omega_\varphi^n=e^{\ee\varphi}\omega_\varphi^k\wedge\Theta_t,
 \qquad\Theta_t=(1-t)\omega^m+t\chi^m,
\end{equation}
put
\begin{equation}\label{eq:C0-normalized-density}
 \nu_t=\omega^k\wedge\Theta_t,\qquad
 h_t=\frac{\nu_t}{\omega^n}
 =(1-t)+t\frac{\omega^k\wedge\chi^m}{\omega^n},\qquad
 \int_M\nu_t=1.
\end{equation}
Since $\omega$ and $\chi$ are uniformly equivalent Poincar\'e-type
metrics, the density $h_t$ has uniform positive upper and lower bounds.
The sole cone and asymptotic hypotheses are the \emph{endpoint}
conditions (B) and (C) of Section~\ref{sec:intro}. In particular,
choose $\delta_*>0$ so that
$1+\delta_*\le\min\{1+\delta,n/k\}$. Then
\begin{equation}\label{eq:C0-normalized-cone}
 n\omega^{n-1}\ge (1+\delta_*)k\omega^{k-1}\wedge\Theta_t
 \qquad (0\le t\le1).
\end{equation}
Moreover, by (C),
\begin{equation}\label{eq:C0-decay}
 |h_t-1|
 =t\left|\frac{\omega^k\wedge\chi^m}{\omega^n}-1\right|
 \le C\rho^{-\eta_0},\qquad 0\le t\le1.
\end{equation}

\begin{proposition}[Uniform $L^\infty$ estimate for the mixed path]
\label{prop:C0}
Assume the endpoint conditions \emph{(B)} and \emph{(C)}.
There exist $\delta_0,C_0>0$, depending only on the fixed
Poincar\'e geometry and the constants in these conditions, such
that every Poincar\'e-type solution $\varphi=\varphi_{\ee,t}$ of
\eqref{eq:C0-normalized-path} with
\begin{equation}\label{eq:C0-smallness}
 \ee\varphi\ge-\delta_0
\end{equation}
satisfies
\begin{equation}\label{eq:C0-conclusion}
 \|\varphi\|_{L^\infty(M)}\le C_0.
\end{equation}
The constants are independent of $t\in[0,1]$, $\ee>0$, and the
particular solution.
\end{proposition}

Write $g=\omega_\varphi=\omega+\ddc\varphi$ and
$H_p(\varphi)=\varphi|\varphi|^{p-1}$. We prove first the
coercive energy and $L^2$ estimates used in the proposition.

\begin{lemma}[Coercive energy estimate]
\label{lem:C0-energy}
Assume the endpoint condition \emph{(B)}. If $\varphi$ solves
\eqref{eq:C0-normalized-path} and satisfies
\eqref{eq:C0-smallness} for sufficiently small $\delta_0>0$,
then for every $p\ge1$,
\begin{align}
 &p\int_M|\varphi|^{p-1}|d\varphi|_\omega^2\,\omega^n
 +p\sum_{j=1}^{k-1}\int_M|\varphi|^{p-1}
 \ii\partial\varphi\wedge\bar\partial\varphi
 \wedge g^j\wedge\omega^{k-1-j}\wedge\Theta_t\notag\\
 &\hspace{25mm}\le C\int_M H_p(\varphi)
 (\omega^n-\omega^k\wedge\Theta_t).
 \label{eq:C0-p-energy}
\end{align}
where $C$ is uniform in $p,t,\ee$.
\end{lemma}

\begin{proof}
Set $\omega_s=\omega+s\ddc\varphi=(1-s)\omega+sg$ and
\[
 Q_s=n\omega_s^{n-1}-k\omega_s^{k-1}\wedge\Theta_t.
\]
By the difference-of-powers identity,
\[
 \omega^n-\omega^k\wedge\Theta_t
 =-\ddc\varphi\wedge\int_0^1Q_s\,ds
 +(e^{\ee\varphi}-1)g^k\wedge\Theta_t.
\]
Since $H_p(\varphi)(e^{\ee\varphi}-1)\ge0$, and the $Q_s$
are closed, Gaffney--Stokes gives
\begin{align*}
 I_p&:=\int_M H_p(\varphi)(\omega^n-\omega^k\wedge\Theta_t)\\
 &\ge-\int_M H_p(\varphi)\ddc\varphi\wedge\int_0^1Q_s\,ds\\
 &=p\int_M|\varphi|^{p-1}
 \ii\partial\varphi\wedge\bar\partial\varphi
 \wedge\int_0^1Q_s\,ds.
\end{align*}
For each fixed $p$, Poincar\'e-type growth and the bounded positive
quasi-coordinate derivatives of $\varphi$ give the required $L^1$
integration-by-parts forms.

We next establish the two endpoint estimates for the same
positive $(m,m)$-form $\Theta_t$. By
\eqref{eq:C0-normalized-cone},
\[
 n\omega^{n-1}\ge(1+\delta_*)k\omega^{k-1}\wedge\Theta_t.
\]
At the other endpoint, the equation and
$\ee\varphi\ge-\delta_0$ imply
\begin{equation}\label{eq:C0-endpoint}
 ng^{n-1}\ge e^{-\delta_0}k g^{k-1}\wedge\Theta_t.
\end{equation}
Indeed, fix any complex hyperplane $E$ in a $g$-unitary frame.
For $\zeta=\omega,\chi$, the sum-of-principal-minors expression
for the positive form $\zeta^m$ gives
\[
 \left.\frac{k g^{k-1}\wedge\zeta^m}{n g^{n-1}}\right|_E
 \le\frac{g^k\wedge\zeta^m}{g^n}.
\]
(The left side includes the natural restriction of both
$(n-1,n-1)$-forms to $E$.) Taking the same nonnegative combination
as in $\Theta_t$, we find
\[
 \left.\frac{k g^{k-1}\wedge\Theta_t}{n g^{n-1}}\right|_E
 \le\frac{g^k\wedge\Theta_t}{g^n}
 =e^{-\ee\varphi}\le e^{\delta_0},
\]
which proves \eqref{eq:C0-endpoint} for all $E$.

For the interpolation between the two endpoints, fix an arbitrary
complex hyperplane $E$ and a positive Hermitian form $A$ on $E$.
For $\zeta=\omega,\chi$, the function
\[
 \Psi_\zeta(A)=\log\frac{nA^{n-1}}
 {kA^{k-1}\wedge(\zeta^m|_E)}
\]
is concave: up to a constant, it is the logarithmic Hessian quotient
$\log\sigma_{n-1}-\log\sigma_{k-1}$ on $E$.
Since
\[
 \Psi_t(A):=\log\frac{nA^{n-1}}
 {kA^{k-1}\wedge(\Theta_t|_E)}
 =-\log\big((1-t)e^{-\Psi_\omega(A)}
                  +t e^{-\Psi_\chi(A)}\big),
\]
it too is concave. To see this directly, with
\[
 p_\omega=\frac{(1-t)e^{-\Psi_\omega}}
 {(1-t)e^{-\Psi_\omega}+te^{-\Psi_\chi}},\qquad
 p_\chi=1-p_\omega,
\]
we have for any real Hermitian direction $B$,
\[
 D^2\Psi_t[B,B]
 =p_\omega D^2\Psi_\omega[B,B]
 +p_\chi D^2\Psi_\chi[B,B]
 -p_\omega p_\chi
 \big(D\Psi_\omega[B]-D\Psi_\chi[B]\big)^2\le0.
\]
Concavity and the two endpoint estimates give, for every
$s\in[0,1]$,
\begin{equation}\label{eq:C0-interpolated}
 n\omega_s^{n-1}\ge
 e^{(1-s)\log(1+\delta_*)-s\delta_0}
 k\omega_s^{k-1}\wedge\Theta_t.
\end{equation}
Put $a(s)=e^{(1-s)\log(1+\delta_*)-s\delta_0}$.
Expanding $\omega_s^{k-1}$ and integrating gives
\begin{equation}\label{eq:C0-Q-lower}
 \int_0^1Q_s\,ds
 \ge\sum_{j=0}^{k-1}\gamma_j
 g^j\wedge\omega^{k-1-j}\wedge\Theta_t,
\end{equation}
where
\[
 \gamma_j=k\binom{k-1}{j}
 \int_0^1(a(s)-1)s^j(1-s)^{k-1-j}\,ds.
\]
At $\delta_0=0$, $a(s)>1$ for $s<1$, so every $\gamma_j>0$.
There are finitely many $j$, and the coefficients depend
continuously on $\delta_0$. Taking $\delta_0>0$ sufficiently small
ensures $\gamma_j\ge\gamma_*>0$ uniformly in $j$ and $t$.
Because $\Theta_t\ge b\omega^m$ in the strongly positive sense
for a uniform $b>0$, the $j=0$ term satisfies
\[
 \ii\partial\varphi\wedge\bar\partial\varphi
 \wedge\omega^{k-1}\wedge\Theta_t
 \ge b'|d\varphi|_\omega^2\omega^n.
\]
Combining the bound for $\int Q_s$ with the bound for $I_p$
proves \eqref{eq:C0-p-energy}.
\end{proof}

\begin{lemma}[Centered $L^2$ estimates]
\label{lem:C0-centered}
Under endpoint condition \emph{(B)}, let $\varphi$ solve
\eqref{eq:C0-normalized-path} with \eqref{eq:C0-smallness}.
Write
\begin{equation}\label{eq:C0-means}
 a_{\ee,t}=\int_M\varphi\,\omega^n,
 \qquad b_{\ee,t}=\int_M\varphi\,\nu_t.
\end{equation}
Then
\begin{equation}\label{eq:C0-centered}
 \|\varphi-a_{\ee,t}\|_{L^2(\omega^n)}
 +\|\varphi-b_{\ee,t}\|_{L^2(\omega^n)}\le C.
\end{equation}
\end{lemma}

\begin{proof}
Set $a=a_{\ee,t}$, $b=b_{\ee,t}$ and $h=h_t$.
Because $\int_M\nu_t=\int_M\omega^n=1$,
\[
 \int_M(1-h)\omega^n=0.
\]
By Lemma~\ref{lem:C0-energy} at $p=1$ and the Poincar\'e
inequality (Lemma~\ref{lem:poincare}),
\begin{align*}
 \|\varphi-a\|_{L^2(\omega^n)}^2
 &\le C\int_M|d\varphi|_\omega^2\omega^n\\
 &\le C\int_M\varphi(1-h)\omega^n\\
 &=C\int_M(\varphi-a)(1-h)\omega^n\\
 &\le C\|\varphi-a\|_{L^2(\omega^n)}
          \|1-h\|_{L^2(\omega^n)}.
\end{align*}
The last factor is uniformly bounded since $h_t$ is, so
$\|\varphi-a\|_{L^2}\le C$. Furthermore,
\[
 |a-b|=\left|\int_M(a-\varphi)h\omega^n\right|
 \le\|\varphi-a\|_{L^2}\|h\|_{L^2}\le C.
\]
Now $\varphi-b=(\varphi-a)+(a-b)$ gives the second estimate.

\end{proof}

\begin{lemma}[Uniform $L^2$ estimate]
\label{lem:C0-L2}
Under endpoint condition \emph{(B)}, solutions of
\eqref{eq:C0-normalized-path} satisfying
\eqref{eq:C0-smallness} obey
\begin{equation}\label{eq:C0-L2}
 \|\varphi\|_{L^2(M,\omega^n)}\le C.
\end{equation}
\end{lemma}

\begin{proof}
As $\Theta_t$ is closed, Gaffney--Stokes and the normalization
imply
\[
 \int_Mg^n=\int_Mg^k\wedge\Theta_t=1.
\]
Apply Jensen's inequality to the probability measures
$g^k\wedge\Theta_t$ and $g^n=e^{\ee\varphi}g^k\wedge\Theta_t$:
\begin{equation}\label{eq:C0-jensen}
 \int_M\varphi\,g^k\wedge\Theta_t\le0,
 \qquad\int_M\varphi\,g^n\ge0.
\end{equation}
The difference-of-powers identity and integration by parts give
\begin{equation}\label{eq:C0-b-identity}
 \int_M\varphi\,g^k\wedge\Theta_t
 =b-\sum_{j=0}^{k-1}\int_M
 \ii\partial\varphi\wedge\bar\partial\varphi
 \wedge g^j\wedge\omega^{k-1-j}\wedge\Theta_t.
\end{equation}
Every summand is nonnegative. By the centered estimate and mass
identity,
\[
 0\le I_1=\int_M(\varphi-a)(1-h_t)\omega^n\le C.
\]
The $p=1$ estimate in Lemma~\ref{lem:C0-energy} controls the
sum in \eqref{eq:C0-b-identity}: its $j\ge1$ terms occur explicitly,
and the $j=0$ term is controlled by the first gradient term,
using $\Theta_t\le C\omega^m$. Thus the first inequality of
\eqref{eq:C0-jensen} implies $b\le C$.
Similarly, integration by parts gives
\[
 \int_M\varphi\,g^n
 =a-\sum_{j=0}^{n-1}\int_M
 \ii\partial\varphi\wedge\bar\partial\varphi
 \wedge g^j\wedge\omega^{n-1-j}\le a.
\]
The second inequality in \eqref{eq:C0-jensen} yields $a\ge0$.
Since Lemma~\ref{lem:C0-centered} gives $|a-b|\le C$, we obtain
$0\le a\le C$. The centered estimate then proves
\eqref{eq:C0-L2}.
\end{proof}

\begin{proof}[Proof of Proposition~\ref{prop:C0}]
Choose
\begin{equation}\label{eq:C0-eta}
 0<\eta<\min\left\{\eta_0,1,\frac1{n-1}\right\}.
\end{equation}
Set $d\mu=\rho^{-\eta}\omega^n$. By the Poincar\'e cusp volume
asymptotics, $\mu(M)<\infty$ and
$\int_M\rho^\eta\omega^n<\infty$. From
\eqref{eq:C0-decay}, $|1-h_t|\le C\rho^{-\eta}$.
For $r\ge2$ set $v=|\varphi|^{r/2}$ and $w=\rho^{-1/2}v$.
The bound $|d\rho|_\omega/\rho\le C$ and the Poincar\'e
inequality yield
\begin{equation}\label{eq:C0-weighted-gradient}
 \int_M\rho(|dw|_\omega^2+|w|^2)\omega^n
 \le C\int_M|dv|_\omega^2\omega^n
 +C\left(\int_Mv\omega^n\right)^2.
\end{equation}
By Lemma~\ref{lem:C0-energy} with $p=r-1$,
\begin{align}
 \int_M|dv|_\omega^2\omega^n
 &=\frac{r^2}{4}\int_M|\varphi|^{r-2}|d\varphi|_\omega^2\omega^n\notag\\
 &\le Cr\int_M|\varphi|^{r-1}|1-h_t|\omega^n
 \le Cr\int_M|\varphi|^{r-1}\,d\mu.
 \label{eq:C0-grad-bound}
\end{align}
Cauchy--Schwarz gives
\begin{equation}\label{eq:C0-mean-term}
 \left(\int_Mv\omega^n\right)^2
 \le\left(\int_M|\varphi|^r\,d\mu\right)
     \left(\int_M\rho^\eta\omega^n\right).
\end{equation}
Now apply the weighted Sobolev inequality
(Lemma~\ref{lem:weighted-sob}) to $w$ with $q=2$,
$p=2(1+\eta)$. Its hypothesis is satisfied because
$\eta<1/(n-1)$. We obtain
\begin{align}
 &\left(\int_M|\varphi|^{r(1+\eta)}\,d\mu\right)^{1/(1+\eta)}
 \notag\\[-2pt]
 &\qquad\le C\int_M\rho(|dw|_\omega^2+|w|^2)\omega^n
 \le Cr\int_M\big(|\varphi|^{r-1}+|\varphi|^r\big)\,d\mu.
 \label{eq:C0-moser-step}
\end{align}
Set $A_r=\max\{1,\|\varphi\|_{L^r(d\mu)}\}$.
Since $\mu(M)<\infty$, the preceding estimate implies
\[
 A_{r(1+\eta)}\le(Cr)^{1/r}A_r,\qquad r\ge2,
\]
after enlarging $C$ if necessary. For
$r_j=2(1+\eta)^j$, the infinite product
$\prod_{j\ge0}(Cr_j)^{1/r_j}$ converges.
By Lemma~\ref{lem:C0-L2} and $d\mu\le\omega^n$,
$A_2\le C$. Passing to the limit $j\to\infty$ yields
\[
 \|\varphi\|_{L^\infty(d\mu)}\le C.
\]
The density $\rho^{-\eta}$ is everywhere strictly positive,
so this equals the ordinary essential supremum on $M$.
This proves \eqref{eq:C0-conclusion}.
\end{proof}

\section{Closedness: second-order and higher-order estimates}
\label{sec:high-order}
Throughout this section $1\le k<n$ and $m=n-k$.
We consider the mixed-form path
$\Theta_t=(1-t)\omega^m+t\chi^m$ and its solution
$g=\omega+\ddc\varphi$.  Conditions \emph{(B)} and the $C^0$
estimate of Proposition~\ref{prop:C0} are the inputs for the
second-order estimate.  

We first recall two algebraic inequalities for a
positive K\"ahler form $\eta$.  In an $\eta$-unitary frame in which
$g_{i\bar j}=\lambda_i\delta_{ij}$, put
\begin{equation}\label{eq:C2-operator}
 \mathscr F_\eta(g)=\log\frac{g^n}{g^k\wedge\eta^m}
 =\log\frac{\binom nk\det g}{\sigma_k(g)},\qquad
 F_{\eta,i}=\frac{1}{\lambda_i}
  -\frac{\sigma_{k-1}(\lambda|i)}{\sigma_k(\lambda)}>0.
\end{equation}
The determinant and elementary symmetric functions refer to
$\eta^{-1}g$.  The linearization $L_\eta$ is positive and
homogeneity gives
\begin{equation}\label{eq:C2-euler}
 \sum_i F_{\eta,i}\lambda_i=m,\qquad
 L_\eta\varphi=m-\operatorname{tr}_{F_\eta}\omega.
\end{equation}
In the two lemmas below, abbreviate $\widehat\chi:=\eta$,
$F_i:=F_{\eta,i}$, and $T:=\operatorname{tr}_{\widehat\chi}g$.
The lemmas apply to any such frozen $\eta$. In Proposition~\ref{prop:C2}
we use their matrix forms for the positive mixture $\mathcal F_t$.

\begin{lemma}[Trace-gradient inequality]
\label{lem:C2-gradient}
At a $\widehat\chi$-normal point where $g$ is diagonal,
\begin{equation}\label{eq:C2-gradient}
 \frac1T\sum_aF_a|T_a|^2
 \le\sum_{a,i,j}\frac{F_i}{\lambda_j}
       |g_{i\bar j,a}|^2.
\end{equation}
\end{lemma}
\begin{proof}
The K\"ahler identities give
$T_a=\sum_jg_{j\bar j,a}=\sum_jg_{a\bar j,j}$.
Thus Cauchy--Schwarz with weights $\lambda_j$ gives
\[
 |T_a|^2\le\Big(\sum_j\lambda_j\Big)
                  \sum_j\frac{|g_{a\bar j,j}|^2}{\lambda_j}.
\]
Multiply by $F_a/T$, sum over $a$, and enlarge the sum of
nonnegative terms by allowing every differentiation index.
\end{proof}

\begin{lemma}[Quantitative quotient Hessian inequality]\label{lem:C2-quotient}
Let $g=\operatorname{diag}(\lambda_1,\ldots,\lambda_n)>0$ and
\[
 \mathscr F(g)=\log\det g-\log\sigma_k(g)+\log\binom nk.
\]
For any complex $n\times n$ matrix $B=(B_{i\bar j})$, let
$B^*=\overline B^{\,T}$.  The complex polarization of the
matrix Hessian is defined by
\[
 D_g^2\mathscr F[B,B^*]
 :=\left.\frac{\partial^2}{\partial s\,\partial\bar s}
 \mathscr F(g+sB+\bar s B^*)\right|_{s=0},
 \qquad s\in\mathbb C.
\]
Here differentiation is with respect to the matrix variable $g$;
$g+sB+\bar sB^*$ is positive Hermitian for sufficiently small
$|s|$.  Then
\begin{equation}\label{eq:C2-hessian}
 -D_g^2\mathscr F[B,B^*]
 \ge\sum_{i,j}\frac{F_i}{\lambda_j}|B_{i\bar j}|^2.
\end{equation}

\end{lemma}
\begin{proof}
We give the complete matrix calculation, since the quantitative
coefficient on the right-hand side of \eqref{eq:C2-hessian} is
stronger than ordinary concavity.  Throughout the proof the matrix
$g$ is diagonal, but the complex matrix $B$ need not be Hermitian.

\medskip\noindent
\emph{Step 1. Principal minors and their weights.}
Let
\[
 \mathcal I_k=\{I\subset\{1,\ldots,n\}: |I|=k\},
 \qquad g_I=(g_{i\bar j})_{i,j\in I},
 \qquad B_I=(B_{i\bar j})_{i,j\in I}.
\]
By the principal-minor formula for the elementary symmetric
polynomial,
\[
 \sigma_k(g)=\sum_{I\in\mathcal I_k}\det g_I,
 \qquad \det g_I=\prod_{i\in I}\lambda_i.
\]
Define
\[
 w_I:=\frac{\det g_I}{\sigma_k(g)}
      =\frac{\prod_{i\in I}\lambda_i}{\sigma_k(\lambda)},
 \qquad
 p_i:=\sum_{I\ni i}w_I,
 \qquad
 p_{ij}:=\sum_{I\supset\{i,j\}}w_I,
 \quad p_{ii}:=p_i.
\]
All the $w_I$ are positive and $\sum_Iw_I=1$.  In particular,
$0\le p_{ij}\le p_i\le1$.  Grouping the products according to
whether they contain the index $i$ gives
\[
 p_i=\frac{\lambda_i\sigma_{k-1}(\lambda|i)}{\sigma_k(\lambda)},
 \qquad
 F_i=\frac1{\lambda_i}
       -\frac{\sigma_{k-1}(\lambda|i)}{\sigma_k(\lambda)}
     =\frac{1-p_i}{\lambda_i}.
\]
For $i\ne j$ and $k\ge2$ one also has
\[
 p_{ij}=\frac{\lambda_i\lambda_j
      \sigma_{k-2}(\lambda|i,j)}{\sigma_k(\lambda)};
\]
when $k=1$ and $i\ne j$, the defining sum is empty and
$p_{ij}=0$.

\medskip\noindent
\emph{Step 2. Differentiate a single principal minor twice.}
For $I\in\mathcal I_k$ put
\[
 Z_I:=\operatorname{tr}(g_I^{-1}B_I)
      =\sum_{i\in I}\frac{B_{i\bar i}}{\lambda_i}.
\]
Introduce a complex parameter $s$ and the Hermitian perturbation
\[
 g(s)=g+sB+\bar s B^*, \qquad
 g_I(s)=g_I+sB_I+\bar s B_I^*.
\]
For $|s|$ sufficiently small these matrices are positive definite.
The first derivative of the determinant is
\[
 \left.\partial_s\det g_I(s)\right|_{s=0}
   =\det g_I\,\operatorname{tr}(g_I^{-1}B_I)
   =\det g_I\,Z_I,
\]
and similarly
$\left.\partial_{\bar s}\det g_I(s)\right|_{s=0}
   =\det g_I\,\overline{Z_I}$.
Differentiating once more and using
\[
 \left.\partial_{\bar s}\bigl(g_I(s)^{-1}\bigr)\right|_{s=0}
  =-g_I^{-1}B_I^*g_I^{-1},
\]
we obtain
\begin{align*}
 \left.\partial_s\partial_{\bar s}\det g_I(s)\right|_{s=0}
 &=\det g_I\,\overline{Z_I}Z_I
   -\det g_I\,\operatorname{tr}
      (g_I^{-1}B_I^*g_I^{-1}B_I)\\
 &=\det g_I\left(
     |Z_I|^2
    -\operatorname{tr}(g_I^{-1}B_Ig_I^{-1}B_I^*)
   \right).
\end{align*}
Since $g_I$ is diagonal and
$(B^*)_{j\bar i}=\overline{B_{i\bar j}}$, the trace here is
\[
 \operatorname{tr}(g_I^{-1}B_Ig_I^{-1}B_I^*)
 =\sum_{i,j\in I}\frac{|B_{i\bar j}|^2}{\lambda_i\lambda_j}.
\]
Hence the full second-derivative formula for a principal minor is
\begin{equation}\label{eq:C2-minor-hessian}
 D_g^2(\det g_I)[B,B^*]
 =\det g_I\left(
   |Z_I|^2-\sum_{i,j\in I}
          \frac{|B_{i\bar j}|^2}{\lambda_i\lambda_j}
              \right).
\end{equation}

\medskip\noindent
\emph{Step 3. Differentiate $\sigma_k$ and $\log\sigma_k$.}
Summing the first derivatives of the principal minors yields
\begin{equation}\label{eq:C2-sigma-first}
 \frac{D_g\sigma_k[B]}{\sigma_k(g)}
 =\sum_{I\in\mathcal I_k}w_IZ_I.
\end{equation}
Similarly, summing \eqref{eq:C2-minor-hessian} and interchanging
the finite sums gives
\begin{align*}
 \frac{D_g^2\sigma_k[B,B^*]}{\sigma_k(g)}
 &=\sum_Iw_I|Z_I|^2
   -\sum_Iw_I\sum_{i,j\in I}
                \frac{|B_{i\bar j}|^2}{\lambda_i\lambda_j}\\
 &=\sum_Iw_I|Z_I|^2
   -\sum_{i,j}\left(\sum_{I\supset\{i,j\}}w_I\right)
                \frac{|B_{i\bar j}|^2}{\lambda_i\lambda_j}\\
 &=\sum_Iw_I|Z_I|^2
   -\sum_{i,j}p_{ij}\frac{|B_{i\bar j}|^2}
                             {\lambda_i\lambda_j}.
\end{align*}
Applying the chain rule for the logarithm, and using
\eqref{eq:C2-sigma-first} and its conjugate, we find
\begin{align}
 D_g^2(\log\sigma_k)[B,B^*]
 &=\frac{D_g^2\sigma_k[B,B^*]}{\sigma_k(g)}
   -\frac{D_g\sigma_k[B]D_g\sigma_k[B^*]}{\sigma_k(g)^2}
     \notag\\
 &=\sum_Iw_I|Z_I|^2
   -\left|\sum_Iw_IZ_I\right|^2
   -\sum_{i,j}p_{ij}\frac{|B_{i\bar j}|^2}
                                {\lambda_i\lambda_j}.
 \label{eq:C2-logsigma-hessian}
\end{align}
The first two terms on the last line form a nonnegative weighted
variance.  Indeed, if $\overline Z=\sum_Iw_IZ_I$, then
\begin{equation}\label{eq:C2-weighted-variance}
 \sum_Iw_I|Z_I|^2-\left|\sum_Iw_IZ_I\right|^2
 =\sum_Iw_I|Z_I-\overline Z|^2\ge0,
\end{equation}
where we used $\sum_Iw_I=1$.

\medskip\noindent
\emph{Step 4. Combine with the determinant and compare coefficients.}
The second derivative of the logarithmic determinant is
\[
 D_g^2(\log\det g)[B,B^*]
 =-\operatorname{tr}(g^{-1}Bg^{-1}B^*)
 =-\sum_{i,j}\frac{|B_{i\bar j}|^2}{\lambda_i\lambda_j}.
\]
Since $\log\binom nk$ is constant, subtracting
\eqref{eq:C2-logsigma-hessian} gives the \emph{exact identity}
\begin{align}
 -D_g^2\mathscr F[B,B^*]
 &=\sum_{i,j}\frac{1-p_{ij}}{\lambda_i\lambda_j}
                  |B_{i\bar j}|^2
   +\sum_Iw_I|Z_I|^2
    -\left|\sum_Iw_IZ_I\right|^2\notag\\
 &=\sum_{i,j}\frac{1-p_{ij}}{\lambda_i\lambda_j}
                  |B_{i\bar j}|^2
   +\sum_Iw_I|Z_I-\overline Z|^2.
 \label{eq:C2-hessian-exact}
\end{align}
The second sum is nonnegative by
\eqref{eq:C2-weighted-variance}.  Moreover,
$p_{ij}\le p_i$ and $F_i=(1-p_i)/\lambda_i$ imply
\[
 \frac{1-p_{ij}}{\lambda_i\lambda_j}
 \ge\frac{1-p_i}{\lambda_i\lambda_j}
 =\frac{F_i}{\lambda_j}.
\]
Using these two facts in \eqref{eq:C2-hessian-exact} proves
\eqref{eq:C2-hessian}.  If $B$ is Hermitian, then $B^*=B$;
thus \eqref{eq:C2-hessian} also implies the ordinary concavity of
$\mathscr F$ on the positive Hermitian cone.

\end{proof}

\begin{proposition}[Second-order estimate]
\label{prop:C2}
Suppose the endpoint strict cone condition \emph{(B)} holds:
\[
 n\omega^{n-1}\ge (1+\delta)k\omega^{k-1}\wedge\chi^m,
 \qquad\delta>0.
\]
Fix $C_1>0$.  There exist $\varepsilon_0>0$ and $A,C>0$,
depending only on $n,k$, the fixed Poincar\'e geometry, $\delta$,
and $C_1$, such that, whenever $0<\varepsilon<\varepsilon_0$
and $t\in[0,1]$, every Poincar\'e-type solution of
\[
 g^n=e^{\varepsilon\varphi}g^k\wedge\Theta_t,\qquad
 g=\omega+\ddc\varphi>0,\qquad
 \Theta_t=(1-t)\omega^m+t\chi^m,
\]
with $\|\varphi\|_{L^\infty(M)}\le C_1$ satisfies
\begin{equation}\label{eq:C2-trace}
 \operatorname{tr}_\chi g\le
 C\exp\!\left(A(\varphi-\inf_M\varphi)\right).
\end{equation}
In particular, for a uniform constant $C'>1$,
\begin{equation}\label{eq:C2-equivalence}
 (C')^{-1}\chi\le g\le C'\chi.
\end{equation}
Only the endpoint cone condition \emph{(B)} enters this
second-order argument; condition \emph{(C)} is used in the
preceding $C^0$ estimate.
\end{proposition}
\begin{proof}
Set
\[
 \mathcal F_t(G,x)=\log\frac{G^n}{G^k\wedge\Theta_t(x)},
 \qquad F^{a\bar b}=D_G\mathcal F_t|_{G=g},\qquad
 L=F^{a\bar b}\partial_a\partial_{\bar b}.
\]
Here $D_G$ differentiates the positive Hermitian matrix
variable, while $\partial^x$ differentiates only the fixed
background coefficients with $G$ held constant.  We use the
contraction convention $D_G\mathcal F_t[B]=F^{a\bar b}B_{a\bar b}$.
The identity
\[
 \mathcal F_t=-\log\big((1-t)e^{-\mathscr F_\omega}
                          +te^{-\mathscr F_\chi}\big)
\]
shows that $F$ is positive definite.  Both frozen quotient
operators are homogeneous of degree $m$ before taking logarithms;
consequently
\begin{equation}\label{eq:C2-new-euler}
 F^{a\bar b}g_{a\bar b}=m,\qquad
 L\varphi=m-\operatorname{tr}_F\omega,
 \qquad \operatorname{tr}_F\omega=F^{a\bar b}\omega_{a\bar b}.
\end{equation}
All calculations below take place in $\chi$-normal holomorphic
coordinates, chosen to diagonalize $g$ at the point under
consideration:
\[
 \chi_{a\bar b}=\delta_{ab},\qquad
 \partial_i\chi_{a\bar b}=0,\qquad
 g_{a\bar b}=\lambda_a\delta_{ab}.
\]
\emph{The matrix $F^{a\bar b}$ need not be diagonal} in this frame.
Write $T=\operatorname{tr}_\chi g=\sum_i\lambda_i$ and
$S=\operatorname{tr}_\chi F=\sum_aF^{a\bar a}>0$.

\medskip\noindent
\textit{Step 1. Construct the auxiliary function and compute its
linearized trace.}
For a constant $A>0$ to be determined, let
\begin{equation}\label{eq:C2-H}
 H=\log T-A\varphi.
\end{equation}
By \eqref{eq:C2-new-euler},
\begin{equation}\label{eq:C2-H-basic}
 LH=\frac{LT}{T}
 -\frac{F^{a\bar b}T_aT_{\bar b}}{T^2}
 +A(\operatorname{tr}_F\omega-m).
\end{equation}
In our normal coordinates, with
$R^\chi_{a\bar b i\bar j}
=-\partial_a\partial_{\bar b}\chi_{i\bar j}$, we have
\[
 T_{a\bar b}
 =\sum_i g_{i\bar i,a\bar b}
  +\sum_i\lambda_i R^\chi_{a\bar b i\bar i}.
\]
By the K\"ahler symmetry
$g_{i\bar i,a\bar b}=g_{a\bar b,i\bar i}$,
\begin{equation}\label{eq:C2-LT-start}
 LT=\sum_iF^{a\bar b}g_{a\bar b,i\bar i}
     +\sum_{a,b,i}F^{a\bar b}\lambda_i
          R^\chi_{a\bar b i\bar i}.
\end{equation}
Put $B_i=(g_{a\bar b,i})_{a,b}$.  Twice differentiating
$\mathcal F_t(g,x)=\varepsilon\varphi$ yields the exact chain rule
\begin{align}
 \varepsilon\varphi_{i\bar i}
 ={}&F^{a\bar b}g_{a\bar b,i\bar i}
       +D_g^2\mathcal F_t[B_i,B_i^*]\notag\\
    &+2\operatorname{Re}
           D_g(\partial_{\bar i}^x\mathcal F_t)[B_i]
       +\partial_i^x\partial_{\bar i}^x\mathcal F_t.
 \label{eq:C2-Fsecond}
\end{align}

Summing \eqref{eq:C2-Fsecond} in $i$ and using
$\sum_i\varphi_{i\bar i}=T-\operatorname{tr}_\chi\omega$
provides the fourth-order cancellation:
\begin{align}
 LT={}&\varepsilon(T-\operatorname{tr}_\chi\omega)
   -\sum_i D_g^2\mathcal F_t[B_i,B_i^*]\notag\\
 &-2\operatorname{Re}\sum_i
          D_g(\partial_{\bar i}^x\mathcal F_t)[B_i]
   -\sum_i\partial_i^x\partial_{\bar i}^x\mathcal F_t
   +\mathcal R_\chi,
 \label{eq:C2-LT-exact}
\end{align}
where
$\mathcal R_\chi=\sum_{a,b,i}F^{a\bar b}\lambda_i
 R^\chi_{a\bar b i\bar i}$.
Since $\chi$ has uniformly bounded curvature,
\begin{equation}\label{eq:C2-R-bound}
 |\mathcal R_\chi|\le CT S,
\end{equation}
where $C$ only depends on $n,\chi.$

\medskip\noindent
\textit{Step 2. Estimate $LH$.}
We first establish the algebraic inequality needed to handle the
additional terms in \eqref{eq:C2-Fsecond}.  Let $\eta$ be either
$\omega$ or $\chi$, and define
\begin{equation}\label{eq:C2-new-K}
 \mathcal K_\eta(B)
 :=-D_g^2\mathscr F_\eta[B,B^*]
     -D_g\mathscr F_\eta[B g^{-1}B^*]\ge0.
\end{equation}

Subtracting the quadratic term from the exact identity
\eqref{eq:C2-hessian-exact} gives
\begin{equation}\label{eq:C2-new-K-exact}
 \mathcal K_\eta(B)
 =\sum_{i,j}\frac{p_i-p_{ij}}{\lambda_i\lambda_j}
       |B_{i\bar j}|^2
    +\sum_Iw_I|Z_I-\overline Z|^2\ge0.
\end{equation}
In particular, $\mathcal K_\eta$ is a nonnegative Hermitian
quadratic form on complex matrices.  Write
$\mathcal K_\eta(B,C)$ for its sesquilinear polarization,
so that
\begin{equation}\label{eq:C2-new-K-CS}
 |\mathcal K_\eta(B,C)|^2
 \le\mathcal K_\eta(B)\mathcal K_\eta(C).
\end{equation}
For any matrix $U$, the coefficients in
\eqref{eq:C2-new-K-exact} give the uniform estimate
\begin{equation}\label{eq:C2-new-K-special}
 \mathcal K_\eta(U^*g)\le C(n,k)|U|_\eta^2.
\end{equation}
Indeed, the off-diagonal coefficient for
$(U^*g)_{i\bar j}=\overline{U_{j\bar i}}\lambda_j$ is
\[
 \frac{\lambda_j}{\lambda_i}(p_i-p_{ij})
 =\frac{\lambda_j\sigma_{k-1}(\lambda|i,j)}
             {\sigma_k(\lambda)}\le1 \quad (i\ne j),
\]
while the diagonal coefficients vanish and the weighted variance
of $Z_I=\sum_{i\in I}\overline{U_{i\bar i}}$ is at most
$C(n,k)|U|_\eta^2$.

We next differentiate the simultaneous-congruence invariance
\[
 \mathscr F_{C^*\eta C}(C^*gC)=\mathscr F_\eta(g),
 \qquad C=I+sU,
\]
first in the complex parameter $s$ and then in the matrix
direction $B$.  With
$\mathcal K_\eta(B,C)$ defined above, this yields
\begin{equation}\label{eq:C2-new-gauge}
 D_gD_\eta\mathscr F_\eta[B,\eta U]
 =\mathcal K_\eta(B,U^*g).
\end{equation}
Here derivatives are complex-polarized as in
Lemma~\ref{lem:C2-quotient}; the identity follows directly from
\[
 D_g^2\mathscr F_\eta[B,gU]
 +D_g\mathscr F_\eta[BU]
 +D_gD_\eta\mathscr F_\eta[B,\eta U]=0.
\]
Taking $U=\eta^{-1}\partial_a\eta$ in
\eqref{eq:C2-new-gauge}, and using
\eqref{eq:C2-new-K-CS}--\eqref{eq:C2-new-K-special} and
bounded Poincar\'e background derivatives, gives
\begin{equation}\label{eq:C2-new-mixed-frozen}
 \big|D_g(\partial_a^x\mathscr F_\eta)[B]\big|
 \le C(n,k,\eta)\sqrt{\mathcal K_\eta(B)}.
\end{equation}

Now put $f_0=\mathscr F_\omega$, $f_1=\mathscr F_\chi$ and
\[
 p_0=\frac{(1-t)e^{-f_0}}
            {(1-t)e^{-f_0}+te^{-f_1}},\qquad
 p_1=\frac{te^{-f_1}}
            {(1-t)e^{-f_0}+te^{-f_1}}.
\]
Then $p_j\ge0$, $p_0+p_1=1$ and
\[
 D_g\mathcal F_t=p_0D_g f_0+p_1D_g f_1.
\]
Direct differentiation of the negative log-sum-exp gives
\begin{equation}\label{eq:C2-new-K-mixture}
 \begin{split}
 \mathcal K_t(B)&:=
 -D_g^2\mathcal F_t[B,B^*]
   -D_g\mathcal F_t[B g^{-1}B^*]\\
 &=p_0\mathcal K_\omega(B)+p_1\mathcal K_\chi(B)
   +p_0p_1\big|D_gf_0[B]-D_gf_1[B]\big|^2.
 \end{split}
\end{equation}
For the first spatial derivatives,
$\partial_a^x\mathcal F_t
=p_0\partial_a^x f_0+p_1\partial_a^x f_1$.
The individual spatial derivatives $\partial_a^x f_j$ are
bounded uniformly for all positive $g$: this follows by
writing the derivative of $g^k\wedge\eta^m$ with $g$ frozen
and using the uniform boundedness of $\eta^{-1}\partial\eta$.
Differentiating the weights $p_j$ in $g$, and combining
\eqref{eq:C2-new-mixed-frozen} with the last variance term in
\eqref{eq:C2-new-K-mixture}, yields
\begin{equation}\label{eq:C2-new-mixed}
  \big|D_g(\partial_a^x\mathcal F_t)[B]\big|
 \le C(n,k,\omega,\chi)\sqrt{\mathcal K_t(B)}.
\end{equation}
For clarity, the derivative of the weights contributes exactly
\[
 -p_0p_1(\partial_a^x f_0-\partial_a^x f_1)
         (D_g f_0[B]-D_g f_1[B]),
\]
which is controlled by the final term in
\eqref{eq:C2-new-K-mixture}.  Similarly, the bounded first and
second background derivatives of $\omega,\chi$, together with
positivity of the mixed wedge products, imply
\begin{equation}\label{eq:C2-new-background2}
 |\partial_i^x\partial_{\bar i}^x\mathcal F_t|\le C(n,k,\chi,\omega).
\end{equation}
Notice that the bounds in
\eqref{eq:C2-new-mixed}--\eqref{eq:C2-new-background2}
do not require any a priori upper or lower bound on $g$.

Using $x-2C\sqrt{x}\ge-C^2$ for $x\ge0$,
\eqref{eq:C2-new-K}--\eqref{eq:C2-new-mixed} give
\begin{align}
 &-D_g^2\mathcal F_t[B_i,B_i^*]
   -2\operatorname{Re}
             D_g(\partial_{\bar i}^x\mathcal F_t)[B_i]\notag\\
 &\qquad\ge
 D_g\mathcal F_t[B_i g^{-1}B_i^*]
   +\mathcal K_t(B_i)-2C\sqrt{\mathcal K_t(B_i)}\notag\\
 &\qquad\ge D_g\mathcal F_t[B_i g^{-1}B_i^*]-C.
 \label{eq:C2-new-third-absorb}
\end{align}
On the other hand, the K\"ahler identities yield
$T_a=\sum_jg_{a\bar j,j}$ in our $\chi$-normal coordinates.
Cauchy--Schwarz in the positive Hermitian form $F$ with weights
$\lambda_j$ yields the \emph{matrix version} of
Lemma~\ref{lem:C2-gradient}:
\begin{equation}\label{eq:C2-new-gradient}
  \frac1T F^{a\bar b}T_aT_{\bar b}
 \le\sum_iD_g\mathcal F_t[B_i g^{-1}B_i^*].
\end{equation}
Combining \eqref{eq:C2-LT-exact},
\eqref{eq:C2-R-bound},
\eqref{eq:C2-new-background2},
\eqref{eq:C2-new-third-absorb}, and
\eqref{eq:C2-new-gradient}, we conclude
\begin{equation}\label{eq:C2-LT-bound}
 LT\ge\varepsilon(T-\operatorname{tr}_\chi\omega)
  +\frac{F^{a\bar b}T_aT_{\bar b}}{T}
  -CTS-C.
\end{equation}
Substituting this into \eqref{eq:C2-H-basic} gives
\[
 LH\ge A(\operatorname{tr}_F\omega-m)-CS
  -\frac CT
  +\varepsilon\left(1-\frac{\operatorname{tr}_\chi\omega}{T}\right).
\]
By \eqref{eq:C2-new-euler},
$m=F^{a\bar b}g_{a\bar b}\le TS$, hence $T^{-1}\le S/m$.
Because $\operatorname{tr}_\chi\omega\le C$ and
$0<\varepsilon\le1$, the last two terms are bounded below
by $-CS$.  Thus
\begin{equation}\label{eq:C2-LH}
  LH\ge A(\operatorname{tr}_F\omega-m)-CS,
 \qquad S=\operatorname{tr}_\chi F.
\end{equation}

\medskip\noindent
\textit{Step 3. Use the endpoint cone condition and a rank-one
comparison at a large eigenvalue.}
From \emph{(B)} and $n>k$, choose $\delta_*>0$ with
$1+\delta_*\le\min\{1+\delta,n/k\}$.  Since
$\Theta_t=(1-t)\omega^m+t\chi^m$,
\begin{equation}\label{eq:C2-B-implies-faceform}
  n\omega^{n-1}
 \ge(1+\delta_*)k\omega^{k-1}\wedge\Theta_t
 \quad(0\le t\le1).
\end{equation}
In the chosen $\chi$-unitary frame diagonalizing $g$,
let $E_i$ be the coordinate hyperplane perpendicular to the
$i$th eigenvector.  For a positive $(1,1)$-form $U$ put
\[
 \mathcal Q_{t,i}(U)=
 \frac{n(U|_{E_i})^{n-1}}
      {k(U|_{E_i})^{k-1}\wedge(\Theta_t|_{E_i})}.
\]
Equation \eqref{eq:C2-B-implies-faceform} shows that
$\mathcal Q_{t,i}(\omega)\ge1+\delta_*$.
The map $U\mapsto\mathcal Q_{t,i}(U)$ is monotone in the
positive Hermitian matrix order and homogeneous of degree
$m$.  Indeed it is the positive reciprocal of a convex
combination of the inverse $(n-1)$-dimensional Hessian
quotients with the fixed $(m,m)$-forms
$\omega^m|_{E_i},\chi^m|_{E_i}$; each corresponding logarithmic
quotient is increasing.  

As $\omega\ge a_0\chi$ for some uniform $a_0>0$, choose
$\theta\in(0,a_0)$ so small that
\[
 (1-\theta/a_0)^m(1+\delta_*)
 \ge1+\tfrac34\delta_*.
\]
Set $W=\omega-\theta\chi>0$.  Then
$W\ge(1-\theta/a_0)\omega$, hence
\begin{equation}\label{eq:C2-W-face}
 \mathcal Q_{t,i}(W)\ge1+\tfrac34\delta_*,
 \qquad\text{uniformly in }x,t,i.
\end{equation}
Let $P_i=\ii\,dz_i\wedge d\bar z_i$ in the selected
$\chi$-unitary coordinates and write $Y_{i,R}=W+RP_i$.
The wedge expansions, valid also when $W$ does not commute
with $g$, are
\[
 Y_{i,R}^n=W^n+nR\,P_i\wedge W^{n-1},
 \qquad
 Y_{i,R}^k\wedge\Theta_t
 =W^k\wedge\Theta_t+kR\,P_i\wedge W^{k-1}\wedge\Theta_t.
\]
Consequently
\begin{equation}\label{eq:C2-rank-one-limit}
 e^{\mathcal F_t(Y_{i,R})}
 \longrightarrow\mathcal Q_{t,i}(W)
 \qquad(R\to\infty).
\end{equation}
Uniform equivalence of $W,\omega,\chi$ and positivity of
$\Theta_t$ imply that this convergence is $O(R^{-1})$
uniformly in the point, parameter, and unitary frame.
Choose one fixed $R$ sufficiently large so that
\begin{equation}\label{eq:C2-rank-one-gap}
 \mathcal F_t(Y_{i,R})\ge\log(1+\delta_*/2).
\end{equation}
Reduce $\varepsilon_0$ so that
$e^{\varepsilon_0C_1}\le1+\delta_*/4$.  Since
$\mathcal F_t(g)=\varepsilon\varphi$, we have
$\mathcal F_t(g)\le\log(1+\delta_*/4)$.
Set $\kappa=\log\big((1+\delta_*/2)/(1+\delta_*/4)\big)>0$.
The concavity of $\mathcal F_t$ therefore gives
\begin{align}
 \kappa
 &\le\mathcal F_t(Y_{i,R})-\mathcal F_t(g)
 \le D_g\mathcal F_t[Y_{i,R}-g]\notag\\
 &=\operatorname{tr}_F\omega-m-\theta S+RF^{i\bar i}.
 \label{eq:C2-rank-one-concavity}
\end{align}
For each frozen quotient, its linearization is bounded above
by $g^{-1}$ as a positive Hermitian form.  The same is true of
the convex combination $F$; hence
$0<F^{i\bar i}\le\lambda_i^{-1}$ in our diagonal frame.
Choosing $i$ with $\lambda_i=\lambda_{\max}$ yields
\begin{equation}\label{eq:C2-subsolution-quant}
  \lambda_{\max}\ge\frac{2R}{\kappa}
 \quad\Longrightarrow\quad
 \operatorname{tr}_F\omega-m\ge\theta S+\frac\kappa2.
\end{equation}

\medskip\noindent
\textit{Step 4. Apply the generalized maximum principle and
obtain the two-sided metric estimate.}
For a fixed Poincar\'e-type solution, $T$ and $H$ are
qualitatively bounded, and $L$ has coefficients bounded in
quasi-coordinates.  The generalized maximum principle
(Proposition~\ref{prop:gen-max}) thus produces a sequence
$x_\ell\in M$ such that
\[
 H(x_\ell)\to\sup_MH,\qquad LH(x_\ell)\le o(1).
\]
Choose $A$ so that $A\theta>C$ in \eqref{eq:C2-LH}.
If $\lambda_{\max}(x_\ell)\ge2R/\kappa$, then
\eqref{eq:C2-LH} and \eqref{eq:C2-subsolution-quant}
would imply
\[
 LH(x_\ell)\ge(A\theta-C)S(x_\ell)
 +\frac{A\kappa}{2}\ge\frac{A\kappa}{2}>0,
\]
which is impossible for large $\ell$.  Thus
$\lambda_{\max}(x_\ell)<2R/\kappa$ eventually and
$T(x_\ell)\le2nR/\kappa$.  Taking the limit in $H$
gives
\[
 \sup_M H\le\log(2nR/\kappa)-A\inf_M\varphi,
\]
which proves \eqref{eq:C2-trace}.  Since $|\varphi|\le C_1$,
there is a uniform $K>1$ such that
\begin{equation}\label{eq:C2-all-upper}
 \lambda_i\le T\le K,\qquad 1\le i\le n.
\end{equation}
To obtain a positive lower bound, uniform background equivalence
gives $\Theta_t\ge b_0\chi^m$ for some $b_0>0$.
The equation implies
\[
 \frac{g^k\wedge\Theta_t}{g^n}
 =e^{-\varepsilon\varphi}\ge
 b_0\frac{g^k\wedge\chi^m}{g^n},
\]
so, using $\binom nm=\binom nk$,
\begin{equation}\label{eq:C2-inverse-eigen}
 \sigma_m(\lambda_1^{-1},\ldots,\lambda_n^{-1})
 \le b_0^{-1}\binom nm\,e^{-\varepsilon\varphi}
 \le M_0
\end{equation}
for a uniform constant $M_0$.  Every $m$-fold product of
inverse eigenvalues is nonnegative, and those containing a
fixed index $i$ give
\[
 \sigma_m(\lambda^{-1})\ge
 \lambda_i^{-1}\sigma_{m-1}(\lambda^{-1}|i)
 \ge\lambda_i^{-1}\binom{n-1}{m-1}K^{-(m-1)}.
\]
Hence
\begin{equation}\label{eq:C2-all-lower}
 \lambda_i\ge
 \frac{\binom{n-1}{m-1}}{M_0K^{m-1}}>0.
\end{equation}
Together \eqref{eq:C2-all-upper}--\eqref{eq:C2-all-lower}
give \eqref{eq:C2-equivalence} and finish the proof.
\end{proof}

\begin{theorem}[Higher-order estimates]\label{thm:high-order}
Under the assumptions of Proposition~\ref{prop:C2}, for every
integer $\ell\ge2$ and $0<\alpha<1$ there is a constant
$C_{\ell,\alpha}$, independent of $t\in[0,1]$ and
$0<\ee<\ee_0$, such that
\begin{equation}\label{eq:Cm-estimate}
 \|\varphi_{\ee,t}\|_{C^{\ell,\alpha}(M)}
 \le C_{\ell,\alpha}.
\end{equation}
The norms are the quasi-coordinate norms of
Section~\ref{subsec:spaces}.
\end{theorem}
\begin{proof}
By Proposition~\ref{prop:C2},
$C^{-1}\chi\le\omega_\varphi\le C\chi$ uniformly.
In each fixed quasi-coordinate polydisc, the equation
$\mathscr F(\omega+\ddc\varphi,\chi)=\ee\varphi$ is therefore
uniformly elliptic and concave, with structural constants and
background coefficient bounds independent of $t$ and $\ee$.
The trace bound and $\|\varphi\|_\infty\le C_1$ first yield
uniform local $C^{1,\alpha}$ estimates by the linear equation
$\Delta_{\chi_t}\varphi
=\operatorname{tr}_{\chi_t}\omega_\varphi
 -\operatorname{tr}_{\chi_t}\omega$
and interior elliptic $W^{2,p}$ estimates (take $p>2n$).
Thus the right-hand side $\ee\varphi$ is uniformly H\"older.
The complex Evans--Krylov estimate, followed by linear Schauder
estimates, gives uniform local $C^{2,\alpha}$ bounds, possibly
starting with a smaller H\"older exponent.
Differentiating the equation and applying interior Schauder estimates
iteratively yields uniform $C^{\ell,\alpha}$ bounds for every
$\ell\ge2$ and every $0<\alpha<1$.
The quasi-coordinate charts have uniformly controlled geometry and
fixed interior radii, so taking the supremum over this atlas proves
\eqref{eq:Cm-estimate}.
\end{proof}

\section{Continuation and proof of the main theorem}
\label{sec:main-proof}
Theorem~\ref{thm:main} is obtained by continuing
\eqref{eq:intro-path} in $t$ for small fixed $\varepsilon>0$
and then letting $\varepsilon\downarrow0$.

\begin{lemma}[Normalized continuation]\label{lem:whole-path}
Suppose $1\le k<n$ and that \emph{(B)} and \emph{(C)} hold for
$\omega,\chi$.  There exist $\delta_*>0$ and $\ee_*>0$,
depending only on the fixed data, such that for each
$0<\ee<\ee_*$ and $t\in[0,1]$ the equation
\begin{equation}\label{eq:sec5-path}
 \omega_{\varphi_{\ee,t}}^n
 =e^{\ee\varphi_{\ee,t}}
 \omega_{\varphi_{\ee,t}}^k\wedge\Theta_t,
 \qquad \omega_{\varphi_{\ee,t}}>0,
\end{equation}
admits a Poincar\'e-type solution satisfying
\begin{equation}\label{eq:sec5-small}
 \ee\|\varphi_{\ee,t}\|_{L^\infty(M)}\le\frac{\delta_*}{2},
\end{equation}
and, uniformly in $\ee$ and $t$,
\begin{equation}\label{eq:sec5-uniform}
 \|\varphi_{\ee,t}\|_{L^\infty(M)}\le C_0,
 \qquad
 \|\varphi_{\ee,t}\|_{C^{\ell,\alpha}(M)}\le C_{\ell,\alpha}
 \quad(\ell\ge2,\ 0<\alpha<1).
\end{equation}

\end{lemma}
\begin{proof}
By \emph{(B)} and \emph{(C)}, Proposition~\ref{prop:C0}
applies uniformly along the normalized path: for some
$\delta_E,C_0>0$, every solution with
$\ee\varphi\ge-\delta_E$ satisfies
$\|\varphi\|_\infty\le C_0$.  Fix
$0<\delta_*<\delta_E$, and let $\ee_{\rm H}>0$ be the
threshold in Proposition~\ref{prop:C2} for this $C_0$.
Choose
\[
 0<\ee_*<\min\left\{\ee_{\rm H},
       \frac{\delta_*}{2\max\{C_0,1\}}\right\}.
\]
Thus any solution with $\ee\|\varphi\|_\infty<\delta_*$
actually satisfies the uniform improvement
\begin{equation}\label{eq:sec5-improvement}
 \|\varphi\|_\infty\le C_0,
 \qquad \ee\|\varphi\|_\infty
 \le\ee C_0<\frac{\delta_*}{2}.
\end{equation}

For fixed $0<\ee<\ee_*$, let
\[
 \mathcal T_{\ee,\delta_*}
 =\big\{t\in[0,1]:\ \eqref{eq:sec5-path}
 \text{ has a Poincar\'e-type solution with }
 \ee\|\varphi\|_\infty<\delta_*\big\}.
\]
This set contains $0$, since $c(0)=1$, $\chi_0=\omega$,
and $\varphi=0$ solves the equation.  It is relatively open
by Proposition~\ref{prop:openness} and the strict improvement
\eqref{eq:sec5-improvement}.  It is relatively closed:
for $t_j\to t_\infty$ in the set, the $C^0$ bound and
Proposition~\ref{prop:C2} give uniform metric equivalence,
while Theorem~\ref{thm:high-order} yields uniform
quasi-coordinate higher-order estimates.  A diagonal subsequence
converges smoothly on compact sets to a Poincar\'e-type solution
at $t_\infty$; the uniform estimates pass to the limit and
\eqref{eq:sec5-improvement} ensures
$\ee\|\varphi_{\ee,t_\infty}\|_\infty\le\delta_*/2<\delta_*$.
Consequently $\mathcal T_{\ee,\delta_*}=[0,1]$.
The same a priori estimates give \eqref{eq:sec5-small}
and \eqref{eq:sec5-uniform}.
\end{proof}

\begin{proof}[Proof of Theorem~\ref{thm:main}]
\emph{Existence.} Choose a sequence $\ee_j\downarrow0$ with
$\ee_j<\ee_*$, where $\ee_*$ is given by
Lemma~\ref{lem:whole-path}. Let $\varphi_j=\varphi_{\ee_j,1}$.
The lemma, Propositions~\ref{prop:C0} and \ref{prop:C2}, and
Theorem~\ref{thm:high-order} give
\begin{equation}\label{eq:sec5-bounds}
 \|\varphi_j\|_{C^{\ell,\alpha}(M)}\le C_{\ell,\alpha}
 \quad (\ell\ge2),\qquad
 C^{-1}\chi\le\omega_{\varphi_j}\le C\chi.
\end{equation}
After passing to a diagonal subsequence,
$\varphi_j\to\varphi$ in $C^\infty$ on compact subsets of $M$.
By \eqref{eq:sec5-bounds}, the limit has uniformly bounded
positive-order quasi-coordinate derivatives and its metric is
uniformly equivalent to the Poincar\'e metric. Since $\varphi$
is bounded, we have $\omega_\varphi\in\PM_\Omega$.
Moreover, $\ee_j\varphi_j\to0$ uniformly. Passing to the limit in
\[
 \omega_{\varphi_j}^n=e^{\ee_j\varphi_j}
       \omega_{\varphi_j}^k\wedge\chi^{n-k}
\]
gives \eqref{eq:intro-target}.

\smallskip
\noindent\emph{Uniqueness.}
Suppose $g_a=\omega+\ddc\varphi_a\in\PM_\Omega$, $a=0,1$,
are two solutions to~\eqref{eq:intro-target}. Let
$w=\varphi_1-\varphi_0$ and
$g_s=(1-s)g_0+sg_1$ for $0\le s\le1$.
At an endpoint $g_a$, diagonalize $\chi$ relative to $g_a$,
with positive eigenvalues $\beta_1,\ldots,\beta_n$.
The equation implies
\begin{equation}\label{eq:sec5-sigma}
 \sigma_{n-k}(\beta)=\binom{n}{n-k}.
\end{equation}
In the coordinate hyperplane perpendicular to the $i$th direction,
the ratio of the coefficients of
$c(1)k g_a^{k-1}\wedge\chi^{n-k}$ and $n g_a^{n-1}$ is
\[
 \frac{c(1)\sigma_{n-k}(\beta|i)}{\binom{n}{n-k}}
 =1-\frac{c(1)\beta_i\sigma_{n-k-1}(\beta|i)}
             {\binom{n}{n-k}}<1.
\]
Because $g_a$ and $\chi$ are uniformly equivalent Poincar\'e-type
metrics, this inequality has a uniform strict margin. Thus, for a
common $\delta'>0$,
\begin{equation}\label{eq:sec5-endpoint-cone}
 n g_a^{n-1}\ge(1+\delta')c(1)k g_a^{k-1}\wedge\chi^{n-k},
 \qquad a=0,1.
\end{equation}
For every complex hyperplane, concavity of the logarithmic
Hessian quotient $\log(\sigma_{n-1}/\sigma_{k-1})$ along the
segment $g_s$ preserves the ratio bound in
\eqref{eq:sec5-endpoint-cone}. Consequently,
\begin{equation}\label{eq:sec5-coercive}
 Q_s:=n g_s^{n-1}-c(1)k g_s^{k-1}\wedge\chi^{n-k}
 \ge\delta' c(1)k g_s^{k-1}\wedge\chi^{n-k}
 \ge c\omega^{n-1}
\end{equation}
for some $c>0$ independent of $s$.
Subtracting the two equations gives the identity of top-degree forms
\begin{equation}\label{eq:sec5-difference}
 0=\big(g_1^n-c(1)g_1^k\wedge\chi^{n-k}\big)
  -\big(g_0^n-c(1)g_0^k\wedge\chi^{n-k}\big)
  =\ddc w\wedge\int_0^1Q_s\,ds.
\end{equation}
Multiply by $w$ and integrate over $M$. Poincar\'e-type growth
implies that $w$ has at most log--log growth and bounded
positive-order covariant derivatives. Thus the integration-by-parts
forms, and their exterior derivatives, are $L^1$; by
Theorem~\ref{thm:gaffney},
\begin{equation}\label{eq:sec5-uniqueness}
 0=\int_M w\,\ddc w\wedge\int_0^1Q_s\,ds
 =-\int_M\ii\partial w\wedge\bar\partial w
                    \wedge\int_0^1Q_s\,ds.
\end{equation}
The strict positivity in \eqref{eq:sec5-coercive} implies
$\partial w=0$ everywhere. Since $M$ is connected, $w$ is constant;
therefore $g_1=g_0$.
\end{proof}

\clearpage
\phantomsection

\vskip10pt
Linwei Zhang, School of Mathematical Sciences,
University of Science and Technology of China, No. 96 Jinzhai Road, Hefei,
Anhui Province, 230026, China;  zhanglinwei@mail.ustc.edu.cn.

\end{document}